\documentclass[a4paper,11pt,twoside,english]{article}
\usepackage{babel}
\usepackage{latexsym,amsfonts}
\usepackage{amsthm}
\usepackage{amsmath}
\usepackage{mathtools}
\usepackage{paralist}
\usepackage{todonotes}

\usepackage{lastpage}
\usepackage{fancyhdr}
\fancypagestyle{plain}{

  \fancyhead[L,C]{}
  \fancyhead[R]{\today}
  \fancyfoot[R,L]{}
  \fancyfoot[C]{\thepage \ of \pageref{LastPage}}
  }

\title{On certain Opial-type results in Ces\`aro spaces of vector-valued functions}
\author{Jan-David Hardtke}
\date{}

\providecommand{\sm}{\setminus}
\providecommand{\ssq}{\subseteq}

\providecommand{\N}{\ensuremath{\mathbb{N}}}
\providecommand{\R}{\ensuremath{\mathbb{R}}}
\providecommand{\eps}{\ensuremath{\varepsilon}}

\providecommand{\keywords}[1]{
  {\let\thefootnote=\relax
  \footnote{{\em Keywords}: #1}}
  \addtocounter{footnote}{-1}
  }

\providecommand{\AMS}[1]{
  {\let\thefootnote=\relax
  \footnote{{\em AMS Subject Classification} (2010): #1}}
  \addtocounter{footnote}{-1}
  }
  
\providecommand{\comment}[1]{
  {\let\thefootnote=\relax
  \footnote{#1}}
  \addtocounter{footnote}{-1}
  }

\providecommand{\address}{
  {\sc \noindent Department of Mathematics \\
  Freie Universit\"at Berlin \\
  Arnimallee 6, 14195 Berlin \\
  Germany \\}
  }

\DeclarePairedDelimiter{\set}{\lbrace}{\rbrace}
\DeclarePairedDelimiter{\paren}{\lparen}{\rparen}

\DeclarePairedDelimiter{\abs}{\lvert}{\rvert}
\DeclarePairedDelimiter{\norm}{\lVert}{\rVert}

\theoremstyle{definition}
\newtheorem{definition}{Definition}[section]
\newtheorem*{definition*}{Definition}

\newtheorem*{example*}{Example}

\newtheorem*{remark*}{Remark}
\theoremstyle{plain}

\newtheorem*{lemma*}{Lemma}
\newtheorem{proposition}[definition]{Proposition}
\newtheorem*{proposition*}{Proposition}
\newtheorem{theorem}[definition]{Theorem}
\newtheorem*{theorem*}{Theorem}
\newtheorem{corollary}[definition]{Corollary}
\newtheorem*{corolary*}{Corollary}

\newenvironment{Proof}[1][\proofname]{\begin{proof}[#1] \setlength{\parindent}{0pt}}{\end{proof}}
\newenvironment{Abstract}{\centering\begin{minipage}{0.8\textwidth} \noindent \small {\sc Abstract.}}{\end{minipage}\par}

\usepackage{color}
\definecolor{darkgreen}{rgb}{0,0.5,0}

\numberwithin{equation}{section}
\newtagform{colored}[\color{blue}]{\color{blue}(}{\color{blue})}
\usetagform{colored}

\usepackage[colorlinks,linkcolor=blue,citecolor=red,urlcolor=darkgreen]{hyperref}
\providecommand{\email}{{\it E-mail address:} \href{mailto:hardtke@math.fu-berlin.de}{\tt hardtke@math.fu-berlin.de}}

\usepackage{amsrefs}

\begin{document}

\maketitle

\begin{Abstract}
Given a Banach space $X$, we consider Ces\`aro spaces $\text{Ces}_p(X)$ of $X$-valued functions over the interval $[0,1]$, where $1\leq p<\infty$.
We prove that if $X$ has the Opial/uniform Opial property, then certain analogous properties also hold for $\text{Ces}_p(X)$.
We also prove a result on the Opial/uniform Opial property of Ces\`aro spaces of vector-valued sequences.
\end{Abstract}
\AMS{46E40 46E30 46B20}
\comment{This work, as well as the paper \cite{hardtke}, are parts of the author's PhD thesis, which is available under 
http://www.diss.fu-berlin.de/diss/receive/FUDISS{\_}thesis{\_}000000099968.}

\section{Introduction}\label{sec:intro}
Let us begin by recalling the definitions of Ces\`aro sequence and function spaces. For $1\leq p<\infty$, the Ces\`aro sequence space $\text{ces}_p$ is 
defined as the space of all sequences $a=(a_n)_{n\in \N}$ of real numbers such that 
\begin{equation*}
\norm{a}_{\text{ces}_p}:=\paren*{\sum_{n=1}^{\infty}\paren*{\frac{1}{n}\sum_{i=1}^n\abs{a_i}}^p}^{1/p}<\infty.
\end{equation*}
$\norm{\cdot}_{\text{ces}_p}$ defines a norm on $\text{ces}_p$. Leibowitz \cite{leibowitz} and Jagers \cite{jagers} proved that $\text{ces}_1=\set*{0}$
and $\text{ces}_p$ is separable and reflexive for $1<p<\infty$. In \cite{bennett} it was proved that for any $p\in (1,\infty)$, the space $\text{ces}_p$
is not isomorphic to $\ell^q$ for any $q\in [1,\infty]$.\par
The Ces\`aro function space $\text{Ces}_p$ on $[0,1]$ is defined in an analogous way as the space of all measurable functions $f:[0,1] \rightarrow \R$ such that
\begin{equation*}
\norm{f}_{\text{Ces}_p}:=\paren*{\int_0^1\paren*{\frac{1}{t}\int_0^t\abs{f(s)}\,\text{d}s}^p\,\text{d}t}^{1/p}<\infty,
\end{equation*}
where, as usual, two functions are identified if they agree a.\,e. $\norm{\cdot}_{\text{Ces}_p}$ defines a norm on $\text{Ces}_p$.\par
Here are some basic results on Ces\`aro function spaces:
\begin{enumerate}[(1)]
\item $\text{Ces}_1$ is the weighted Lebesgue space $L_w^1[0,1]$, where $w(t):=\log(1/t)$.
\item $\text{Ces}_p|_{[0,a]}$ is a subspace of $L^p[0,a]$ for every $p\in [1,\infty)$ and every $a\in (0,1)$,
but not for $a=1$.
\item $\text{Ces}_p$ is separable and nonreflexive for every $p\in [1,\infty)$.
\item For $1<p<\infty$ one has $L^p[0,1]\ssq \text{Ces}_p$ and $\norm{f}_{\text{Ces}_p}\leq q\norm{f}_p$ for all $f\in L^p[0,1]$,
where $q$ is the conjugated exponent to $p$.
\end{enumerate}
These and further results are collected in \cite{astashkin}*{Theorem 1}. Also, by \cite{astashkin}*{Theorem 7}, for $p\in (1,\infty)$
the space $\text{Ces}_p$ is not isomorphic to $L^q[0,1]$ for any $q\in [1,\infty]$. For further information on Ces\`aro function spaces 
see \cites{astashkin2, astashkin, astashkin3} and references therein. For more information on Ces\`aro sequence spaces see, for example, 
the introduction of \cite{astashkin} and references therein. Results on more general types of Ces\`aro function spaces, where the space 
$L^p$ appearing implicitly in the definition of $\text{Ces}_p$ is replaced by a more general function space, can be found for example in 
\cite{lesnik, lesnik2}.\par
Now consider a real Banach space $X$. $X$ is said to have the fixed point property (resp. weak fixed point property) if for every 
closed and bounded (resp. weakly compact) convex subset $C\ssq X$, every nonexpansive mapping $F:C \rightarrow C$ has a fixed point 
(where $F$ is called nonexpansive if $\norm{F(x)-F(y)}\leq\norm{x-y}$ for all $x,y\in C$).\par
A bounded, closed, convex subset $C\ssq X$ is said to have normal structure provided that for each subset $B\ssq C$ which contains at least two 
elements there exists a point $x\in B$ such that
\begin{equation*}
\sup_{y\in B}\norm{x-y}<\operatorname{diam}(B),
\end{equation*}
where $\text{diam}(B)$ denotes the diameter of $B$. The space $X$ itself is said to have normal structure if every bounded, closed, convex 
subset of $X$ has normal structure. It is well known that if $C$ is weakly compact and has normal structure, then every nonexpansive 
mapping $F:C \rightarrow C$ has a fixed point (see e.\,g. \cite{goebel}*{Theorem 2.1}), thus spaces with normal structure have 
the weak fixed point property. For example, every compact, convex set has normal structure (see e.\,g. \cite{prus}*{p.119}) and hence all
finite-dimensional spaces possess normal structure. Also, every space which is uniformly convex in every direction has normal structure 
(see e.\,g. \cite{prus}*{Corollary 5.6}). An example of a Banach space which fails the weak fixed point property is $L^1[0,1]$ (see \cite{alspach}).\par
The space $X$ is said to have the Opial property if 
\begin{equation*}
\limsup_{n\to \infty}\norm{x_n}<\limsup_{n\to \infty}\norm{x_n-x}
\end{equation*}
holds for every weakly null sequence $(x_n)_{n\in \N}$ in $X$ and every $x\in X\setminus\set*{0}$ (one could
as well use $\liminf$ instead of $\limsup$ or assume from the beginning that both limits exist).\par
This property was first considered by Opial in \cite{opial} (starting from the Hilbert spaces as canonical example) to provide 
a result on iterative approximations of fixed points of nonexpansive mappings. It is shown in \cite{opial} that
the spaces $\ell^p$ for $1\leq p<\infty$ enjoy the Opial property, whereas $L^p[0,1]$ for $1<p<\infty, p\neq 2$
fails to have it. Note further that every Banach space with the Schur property (i.\,e. weak and norm convergence of 
sequences coincide) trivially has the Opial property. Also, $X$ is said to have the nonstrict Opial property if it fulfils 
the definition of the Opial property with ``$\leq$'' instead of ``$<$'' (\cite{sims2}, in \cite{garcia-falset1} it is called 
weak Opial property).\footnote{Note that one always has $\limsup\norm{x_n}\leq\limsup\norm{x_n-x}+\norm{x}\leq2\limsup\norm{x_n-x}$
if $(x_n)_{n\in \N}$ converges weakly to zero, since the norm is weakly lower semicontinuous. In general, the constant $2$ is the best 
possible. Consider, for example, in the space $c$ of all convergent sequences (with sup-norm) the weak null sequence $(2e_n)_{n\in \N}$
(where $e_n$ is the sequence whose $n$-th entry is $1$ and all other entries are $0$) and $x=(1,1,1,\dots)$.} It is known that every weakly 
compact convex set in a Banach space with the Opial property has normal structure (see e.\,g. \cite{prus}*{Theorem 5.4}) and thus the Opial 
property implies the weak fixed point property.\par
The notion of uniform Opial property was introduced by Prus in \cite{prus2}: $X$ is said to have the uniform Opial property if
for every $c>0$ there is some $r>0$ such that
\begin{equation*}
1+r\leq\liminf_{n\to \infty}\norm{x_n-x}
\end{equation*}
holds for every $x\in X$ with $\norm{x}\geq c$ and every weakly null sequence $(x_n)_{n\in \N}$ in $X$ with 
$\liminf\norm{x_n}\geq 1$. In \cite{prus2} it was proved that a Banach space is reflexive and has the uniform
Opial property if and only if it has the so called property $(L)$ (see \cite{prus2} for the definition), and that 
$X$ has the fixed point property whenever $X^*$ has property $(L)$.\par
A modulus corresponding to the uniform Opial property was defined in \cite{lin0}:
\begin{equation*}
r_X(c):=\inf\set*{\liminf_{n\to \infty}\norm{x_n-x}-1} \ \ \forall c>0,
\end{equation*}
where the infimum is taken over all $x\in X$ with $\norm{x}\geq c$ and all weakly null sequences $(x_n)_{n\in \N}$ in $X$ with 
$\liminf\norm{x_n}\geq 1$ (if $X$ has the Schur property, we agree to set $r_X(c):=1$ for all $c>0$). Then $X$ has the uniform 
Opial property if and only $r_X(c)>0$ for every $c>0$.\par
Here we will use instead the following equivalent formulation of the uniform Opial property (\cite{khamsi}*{Definition 3.1}):
$X$ has the uniform Opial property if and only if for every $\eps>0$ and every $R>0$ there is some $\eta>0$ such that
\begin{equation*}
\eta+\liminf_{n\to \infty}\norm{x_n}\leq\liminf_{n\to \infty}\norm{x_n-x}
\end{equation*}
holds for all $x\in X$ with $\norm{x}\geq\eps$ and every weakly null sequence $(x_n)_{n\in \N}$ in $X$ with $\limsup\norm{x_n}\leq R$.\par
In \cite{hardtke} the author defined a modulus for this formulation in the following way:
\begin{equation*}
\eta_X(\eps,R):=\inf\set*{\liminf_{n\to \infty}\norm{x_n-x}-\liminf_{n\to \infty}\norm{x_n}} \ \ \forall \eps,R>0,
\end{equation*}
where the infimum is taken over all $x\in X$ with $\norm{x}\geq\eps$ and all weakly null sequences $(x_n)_{n\in \N}$ in $X$ 
with $\limsup\norm{x_n}\leq R$. Thus $X$ has the uniform Opial property if and only if $\eta_X(\eps,R)>0$ for all $\eps,R>0$
(see also \cite{hardtke}*{Lemma 1.1} for a more precise connection between the moduli $r_X$ and $\eta_X$).\par
In \cite{hardtke} the author studied Opial properties in infinite $\ell^p$-sums and also some analogous results for 
Lebesgue-Bochner spaces of vector-valued functions (these spaces cannot have the usual Opial property, as even $L^p[0,1]$ 
for $1<p<\infty$, $p\neq 2$ does not enjoy this property, but if one replaces weak convergence by pointwise weak convergenvce 
(almost everywhere), then some ``Opial-like'' results for Lebesgue-Bochner spaces can be established, see \cite{hardtke} for the 
detailed formulations and proofs).\par
The purpose of this paper is to prove some results analogous to those of \cite{hardtke} for Ces\`aro spaces of vector-valued functions
and Ces\`aro sums. We will start with the latter.

\section{Opial properties in Ces\`aro sums}\label{sec:cesarosums}
Given a sequence $(X_n)_{n\in \N}$ of Banach spaces and $p\in (1,\infty)$, we define the $p$-Ces\`aro sum $\bigl[\bigoplus_{n\in \N}X_n\bigr]_{\text{ces}_p}$ 
of $(X_n)_{n\in \N}$ as the space of all sequences $x=(x_n)_{n\in \N}$ with $x_n\in X_n$ for each $n$ such that $(\norm{x_n})_{n\in \N}\in \text{ces}_p$, equipped 
with the norm
\begin{equation*}
\norm{x}_{\text{ces}_p}:=\norm*{(\norm{x_n})_{n\in \N}}_{\text{ces}_p}=\paren*{\sum_{n=1}^{\infty}\paren*{\frac{1}{n}\sum_{i=1}^n\norm{x_i}}^p}^{1/p}.
\end{equation*}
\indent In \cite{cui} it was proved that $\text{ces}_p$ has the uniform Opial property for every $p\in (1,\infty)$. In \cite{saejung2}*{Theorem 1} Saejung proved 
that $\text{ces}_p$ can be regarded as a subspace of the $\ell^p$-sum $X_p:=\bigl[\bigoplus_{n\in \N}\ell^1(n)\bigr]_p$ (where $\ell^1(n)$ denotes the $n$-dimensional
space with $\ell^1$-norm) via the isometric embedding $T:\text{ces}_p \rightarrow X_p$ defined by
\begin{equation*}
(Ta)(n):=\frac{1}{n}(a_1,\dots,a_n) \ \ \forall n\in \N, \forall a\in \text{ces}_p.
\end{equation*}
In \cite{saejung2}*{Theorem 7} it is proved that the $\ell^p$-sum of any sequence of finite-dimensional spaces has the uniform Opial property
(see also \cite{prus}*{Example 4.23 (2.)} and Corollary 3.14 in \cite{hardtke}). Thus Saejung obtains a new proof that $\text{ces}_p$ has the 
uniform Opial property (\cite{saejung2}*{Corollary 9}).\par
Saejung's embedding idea directly generalises to $\text{ces}_p$-sums. For a given sequence $(X_n)_{n\in \N}$ of Banach spaces we consider the mapping
$S$ from the Ces\`aro sum $\bigl[\bigoplus_{n\in \N}X_n\bigr]_{\text{ces}_p}$ to $\bigl[\bigoplus_{n\in \N}(X_1\oplus_1\dots\oplus_1 X_n)\bigr]_p$ defined by
\begin{equation*}
(Sx)(n):=\frac{1}{n}(x_1,\dots,x_n) \ \ \forall n\in \N, \forall x\in \Bigl[\bigoplus_{n\in \N}X_n\Bigr]_{\text{ces}_p}.
\end{equation*}
Then $S$ is an isometric embedding.\par
In \cite{hardtke}*{Proposition 3.11} the author proved that for any $1\leq p<\infty$ the $\ell^p$-sum of any family of Banach spaces with 
the Opial property/nonstrict Opial property has again the Opial property/nonstrict Opial property. Thus, via the above embedding, we obtain the following result.
\begin{proposition}\label{prop:cesarosumsopial}
Let $p\in (1,\infty)$. If $(X_n)_{n\in \N}$ is a sequence of Banach spaces such that each $X_n$ has the Opial property (nonstrict Opial property),
then $\bigl[\bigoplus_{n\in \N}X_n\bigr]_{\mathrm{ces}_p}$ also has the Opial property (nonstrict Opial property).
\end{proposition}
It was also proved in \cite{hardtke}*{Theorem 3.13} that, for any family $(X_i)_{i\in I}$ of Banach spaces and every $1\leq p<\infty$, the sum $\bigl[\bigoplus_{i\in I}X_i\bigr]_p$
has the uniform Opial property if 
\begin{equation*}
\inf\set*{\eta_{X_J}(\eps,R):J\ssq I \ \text{finite}}>0 \ \ \forall \eps,R>0,
\end{equation*}
where $X_J:=\bigl[\bigoplus_{i\in J}X_i\bigr]_p$ for $J\ssq I$.\par
Using this together with the above embedding, we obtain the following result.
\begin{proposition}\label{prop:cesarosumsunifopial}
Let $p\in (1,\infty)$ and $(X_n)_{n\in \N}$ be a sequence of Banach spaces. Put $Y_m:=\bigl[\bigoplus_{n=1}^m(X_1\oplus_1\dots\oplus_1X_n)\bigr]_p$
for each $m\in \N$. If
\begin{equation*}
\inf_{m\in \N}\eta_{Y_m}(\eps,R)>0 \ \ \forall \eps,R>0,
\end{equation*}
then $\bigl[\bigoplus_{n\in \N}X_n\bigr]_{\mathrm{ces}_p}$ has the uniform Opial property.
\end{proposition}
Note that this implies in particular that $\bigl[\bigoplus_{n\in \N}X_n\bigr]_{\text{ces}_p}$ has the uniform Opial property if each $X_n$ has the Schur property.

\section{Opial-type properties in Ces\`aro spaces of \texorpdfstring{\\}{} vector-valued functions}\label{sec:cesarofunctions} 
Now we consider Ces\`aro spaces of vector-valued functions on $[0,1]$. Usually, for a given Banach space $X$ and a K\"othe function space $E$ (see for instance 
\cite{lin} for the definition), one considers the K\"othe-Bochner space $E(X)$ of all (equivalence classes of) $X$-valued Bochner-measurable functions $f$ such that $\norm{f(\cdot)}\in E$, 
endowed with the norm $\norm{f}_{E(X)}:=\norm{\norm{f(\cdot)}}_E$ (this includes the Lebesgue-Bochner spaces for $E=L^p$). Such spaces have been intensively
studied (see for example the collection of results in \cite{lin}). However, Ces\`aro function spaces are not K\"othe spaces in the usual sense,
since they are not contained in $L^1$ (see point (2) in the introduction). But the Ces\`aro spaces, like K\"othe spaces, satisfy the important 
monotonicity property: if $f\in \text{Ces}_p$ and $g:[0,1] \rightarrow \R$ is measurable with $\abs*{g(t)}\leq\abs*{f(t)}$ a.\,e., then $g\in \text{Ces}_p$ and 
$\norm{g}_{\text{Ces}_p}\leq\norm{f}_{\text{Ces}_p}$.\par
Therefore, given a Banach space $X$, we can still define the space $\text{Ces}_p(X)$ of all (equivalence classes of) Bochner-measurable functions 
$f:[0,1] \rightarrow X$ such that $\norm{f(\cdot)}\in \text{Ces}_p$, equipped with the norm $\norm{f}_{\text{Ces}_p(X)}:=\norm*{\norm{f(\cdot)}}_{\text{Ces}_p}$.\par
We will now prove a result for sequences of functions in $\text{Ces}_p(X)$ which are pointwise a.\,e. convergent to zero with respect to the weak topology of $X$,
where $X$ is assumed to have the nonstrict Opial property. The result is similar to the one obtained in \cite{hardtke}*{Proposition 4.1} for Lebesgue-Bochner spaces. 
The proof also makes use of similar techniques.

\begin{theorem}\label{thm:opialcesaro}
Let $1\leq p<\infty$ and let $X$ be a Banach space with the nonstrict Opial property. Let $(f_n)_{n\in \N}$ be a bounded sequence in $\mathrm{Ces}_p(X)$
such that $(f_n(t))_{n\in \N}$ converges weakly to zero for almost every $t\in [0,1]$. Suppose further that there exists a $g\in \mathrm{Ces}_p$ such
that $\norm{f_n(t)}\to g(t)$ a.\,e. Let $f\in \mathrm{Ces}_p(X)$ and $\varphi(t):=\liminf_{n\to \infty}\norm{f_n(t)-f(t)}$ for $t\in [0,1]$. Then 
\begin{align}\label{eq:opialcesaro1}
&2^{p-1}\int_0^1\frac{1}{t^p}\paren*{\paren*{\int_0^t\varphi(s)\,\text{d}s}^p-\paren*{\int_0^t g(s)\,\text{d}s}^p}\,\text{d}t \nonumber \\
&\leq2^{p-1}\limsup_{n\to \infty}\norm{f_n-f}_{\mathrm{Ces}_p(X)}^p-\limsup_{n\to \infty}\norm{f_n}_{\mathrm{Ces}_p(X)}^p.
\end{align}
In particular, 
\begin{equation}\label{eq:opialcesaro2}
\limsup_{n\to \infty}\norm{f_n}_{\mathrm{Ces}_p(X)}\leq 2^{1-1/p}\limsup_{n\to \infty}\norm{f_n-f}_{\mathrm{Ces}_p(X)} \ \ \forall f\in \mathrm{Ces}_p(X).
\end{equation}
\end{theorem}

\begin{Proof}
Using the identification of $\text{Ces}_1$ with $L_w^1[0,1]$ from \cite{astashkin}*{Theorem 1} (where $w(t)=\log(1/t)$) the assertion for $p=1$ easily follows
from \cite{hardtke}*{Proposition 4.1}. We will therefore assume $p>1$.\par
So let $f\in \text{Ces}_p(X)$. Without loss of generality, we may assume that $\lim_{n\to \infty}\norm{f_n}_{\text{Ces}_p(X)}$
and $\lim_{n\to \infty}\norm{f_n-f}_{\text{Ces}_p(X)}$ exist and also that $\norm{f_n(t)}\to g(t)$ and $f_n(t)\to 0$ weakly for 
every $t\in [0,1]$.\par
Since $M:=\sup_{n\in \N}\norm{f_n}_{\text{Ces}_p(X)}<\infty$ it follows from Fatou's Lemma that $\varphi\in \text{Ces}_p$.\par
Let us put
\begin{equation*}
a:=\norm{\varphi}_{\text{Ces}_p}^p-\norm{g}_{\text{Ces}_p}^p=\int_0^1\frac{1}{t^p}\paren*{\paren*{\int_0^t\varphi(s)\,\text{d}s}^p
-\paren*{\int_0^t g(s)\,\text{d}s}^p}\,\text{d}t.
\end{equation*}
Since $X$ has the nonstrict Opial property we have $\varphi(t)\geq g(t)$ for every $t\in [0,1]$. Hence $a\geq 0$.\par
Let $0<\eps<1$. Denote by $\lambda$ the Lebesgue measure on $[0,1]$.\par
The equi-integrability of finite subsets of $L^1$ enables us to find a $0<\tau<\eps$ such that for every measurable 
set $A\ssq [0,1]$ one has
\begin{equation}\label{opialcesaro eq:A1}
\lambda(A)\leq\tau \ \Rightarrow \ \int_A\paren*{\frac{1}{t}\int_0^tF(s)\,\text{d}s}^p\,\text{d}t\leq\eps \ \ \forall F\in \set*{\norm{f(\cdot)},\varphi,g}.
\end{equation}
Next we choose $0<\theta<\tau$ such that
\begin{equation}\label{opialcesaro eq:A2}
\frac{p\theta}{1-p}\paren*{\paren*{1-\frac{\tau}{3}}^{1-p}-\paren*{\frac{\tau}{3}}^{1-p}}\paren*{\int_0^{1-\frac{\tau}{3}}F(s)\,\text{d}s}^{p-1}\leq\eps
\end{equation}
for $F\in \set*{\varphi,g}$ and then, again by equi-integrability, we find $\delta>0$ such that for every measurable subset $D\ssq [0,1-\frac{\tau}{3}]$ one has
\begin{equation}\label{opialcesaro eq:A3}
\lambda(D)\leq\delta \ \Rightarrow \ \int_DF(s)\,\text{d}s\leq\theta \ \ \forall F\in \set*{\norm{f(\cdot)},\varphi,g}
\end{equation}
(remember that $\text{Ces}_p|_{[0,b]}\ssq L^1[0,b]$ for every $b\in (0,1)$).\par
Now we apply Egorov's theorem (cf. \cite{halmos}*{Theorem A, p.88}) to find a measurable set $C\ssq [0,1]$ with $\lambda([0,1]\sm C)\leq\delta$ such that 
$\norm{f_n(t)}\to g(t)$ uniformly in $t\in C$. It follows that
\begin{equation*}
\lim_{n\to \infty}\int_{[0,t]\cap C}\norm{f_n(s)}\,\text{d}s=\int_{[0,t]\cap C}g(s)\,\text{d}s \ \ \forall t\in [0,1).
\end{equation*}
Thus we can apply Egorov's theorem once more to deduce that there exists a measurable set $F\ssq [0,1)$ with $\lambda([0,1]\sm F)\leq\tau/3$ such that
\begin{equation*}
\lim_{n\to \infty}\paren*{\frac{1}{t}\int_{[0,t]\cap C}\norm{f_n(s)}\,\text{d}s}^p=\paren*{\frac{1}{t}\int_{[0,t]\cap C}g(s)\,\text{d}s}^p \ \ \text{uniformly\ in}\ t\in F.
\end{equation*}
Put $B:=F\cap [\frac{\tau}{3},1-\frac{\tau}{3}]$. Then $\lambda([0,1]\sm B)\leq\lambda([0,1]\sm F)+2\tau/3\leq\tau$ and
\begin{equation}\label{opialcesaro eq:A4}
\lim_{n\to \infty}\int_B\paren*{\frac{1}{t}\int_{[0,t]\cap C}\norm{f_n(s)}\,\text{d}s}^p\,\text{d}t=\int_B\paren*{\frac{1}{t}\int_{[0,t]\cap C}g(s)\,\text{d}s}^p\,\text{d}t.
\end{equation}
Since $M=\sup_{n\in \N}\norm{f_n}_{\text{Ces}_p(X)}<\infty$, we can find a subsequence $(n_k)_{k\in \N}$ of indices such that all the limits involved in the 
following calculations exist.\par
We have
\begin{align*}
&\lim_{n\to \infty}\norm{f_n}_{\text{Ces}_p(X)}^p \\
&=\lim_{k\to \infty}\paren*{\int_B\paren*{\frac{1}{t}\int_0^t\norm{f_{n_k}(s)}\,\text{d}s}^p\,\text{d}t+
\int_{[0,1]\sm B}\paren*{\frac{1}{t}\int_0^t\norm{f_{n_k}(s)}\,\text{d}s}^p\,\text{d}t} \\
&\leq2^{p-1}\lim_{k\to \infty}\paren*{\int_B\paren*{\frac{1}{t}\int_{[0,t]\cap C}\norm{f_{n_k}(s)}\,\text{d}s}^p\,\text{d}t+
\int_B\paren*{\frac{1}{t}\int_{[0,t]\sm C}\norm{f_{n_k}(s)}\,\text{d}s}^p\,\text{d}t} \\
&+\lim_{k\to \infty}\int_{[0,1]\sm B}\paren*{\frac{1}{t}\int_0^t\norm{f_{n_k}(s)}\,\text{d}s}^p\,\text{d}t,
\end{align*}
where we have used the inequality $(a+b)^p\leq2^{p-1}(a^p+b^p)$ for $a,b\geq 0$, which is due to the convexity of the function $t\mapsto t^p$.\par
From \eqref{opialcesaro eq:A4} it now follows that
\begin{align}\label{opialcesaro eq:A5}
&\nonumber \lim_{n\to \infty}\norm{f_n}_{\text{Ces}_p(X)}^p \\
&\nonumber \leq2^{p-1}\paren*{\int_B\paren*{\frac{1}{t}\int_{[0,t]\cap C}g(s)\,\text{d}s}^p\,\text{d}t+ 
\lim_{k\to \infty}\int_B\paren*{\frac{1}{t}\int_{[0,t]\sm C}\norm{f_{n_k}(s)}\,\text{d}s}^p\,\text{d}t} \\
&+\lim_{k\to \infty}\int_{[0,1]\sm B}\paren*{\frac{1}{t}\int_0^t\norm{f_{n_k}(s)}\,\text{d}s}^p\,\text{d}t.
\end{align}
Because of $\lambda([0,1]\sm B)\leq\tau$ and \eqref{opialcesaro eq:A1} we have
\begin{equation*}
\int_{[0,1]\sm B}\paren*{\frac{1}{t}\int_0^t\norm{f(s)}\,\text{d}s}^p\,\text{d}t\leq\eps.
\end{equation*}
Thus by the triangle inequality for $L^p$ we get
\begin{align*}
&\lim_{k\to \infty}\paren*{\int_{[0,1]\sm B}\paren*{\frac{1}{t}\int_0^t\norm{f_{n_k}(s)}\,\text{d}s}^p\,\text{d}t}^{1/p} \\
&\leq\lim_{k\to \infty}\paren*{\int_{[0,1]\sm B}\paren*{\frac{1}{t}\int_0^t\norm{f_{n_k}(s)-f(s)}\,\text{d}s}^p\,\text{d}t}^{1/p}+\eps^{1/p}.
\end{align*}
It follows that
\begin{align*}
&\lim_{k\to \infty}\int_{[0,1]\sm B}\paren*{\frac{1}{t}\int_0^t\norm{f_{n_k}(s)}\,\text{d}s}^p\,\text{d}t \\
&\leq\lim_{k\to \infty}\int_{[0,1]\sm B}\paren*{\frac{1}{t}\int_0^t\norm{f_{n_k}(s)-f(s)}\,\text{d}s}^p\,\text{d}t \\
&+\lim_{k\to \infty}\Bigg|\paren*{\paren*{\int_{[0,1]\sm B}\paren*{\frac{1}{t}\int_0^t\norm{f_{n_k}(s)-f(s)}\,\text{d}s}^p\,\text{d}t}^{1/p}+\eps^{1/p}}^p \\
&-\int_{[0,1]\sm B}\paren*{\frac{1}{t}\int_0^t\norm{f_{n_k}(s)-f(s)}\,\text{d}s}^p\,\text{d}t\Bigg|.
\end{align*}
Put $L:=M+\norm{f}_{\text{Ces}_p(X)}+1$. Since $\abs*{a^p-b^p}\leq pL^{p-1}\abs*{a-b}$ for all $a,b\in [0,L]$ (mean-value theorem) we obtain
\begin{align}\label{opialcesaro eq:A6}
&\nonumber \lim_{k\to \infty}\int_{[0,1]\sm B}\paren*{\frac{1}{t}\int_0^t\norm{f_{n_k}(s)}\,\text{d}s}^p\,\text{d}t \\
&\leq\lim_{k\to \infty}\int_{[0,1]\sm B}\paren*{\frac{1}{t}\int_0^t\norm{f_{n_k}(s)-f(s)}\,\text{d}s}^p\,\text{d}t+pL^{p-1}\eps^{1/p}.
\end{align}
Next we define $h(s):=(1-p)^{-1}(s^p-s/3^{1-p})$ for $s\geq 0$.\par
Recall that $B\ssq [\frac{\tau}{3},1-\frac{\tau}{3}]$, $\lambda([0,1]\sm C)\leq\delta$ and $\theta<\tau$. Thus it follows from \eqref{opialcesaro eq:A3} that
\begin{equation*}
\int_B\paren*{\frac{1}{t}\int_{[0,t]\sm C}\norm{f(s)}\,\text{d}s}^p\,\text{d}t
\leq\tau^p\int_B\frac{1}{t^p}\,\text{d}t\leq\tau^p\int_{\frac{\tau}{3}}^1\frac{1}{t^p}\,\text{d}t=h(\tau).
\end{equation*}
Hence
\begin{align*}
&\lim_{k\to \infty}\paren*{\int_B\paren*{\frac{1}{t}\int_{[0,t]\sm C}\norm{f_{n_k}(s)}\,\text{d}s}^p\,\text{d}t}^{1/p} \\
&\leq\lim_{k\to \infty}\paren*{\int_B\paren*{\frac{1}{t}\int_{[0,t]\sm C}\norm{f_{n_k}(s)-f(s)}\,\text{d}s}^p\,\text{d}t}^{1/p}+h(\tau)^{1/p}.
\end{align*}
Using the same trick as before we now obtain
\begin{align}\label{opialcesaro eq:A7}
&\nonumber \lim_{k\to \infty}\int_B\paren*{\frac{1}{t}\int_{[0,t]\sm C}\norm{f_{n_k}(s)}\,\text{d}s}^p\,\text{d}t \\
&\leq\lim_{k\to \infty}\int_B\paren*{\frac{1}{t}\int_{[0,t]\sm C}\norm{f_{n_k}(s)-f(s)}\,\text{d}s}^p\,\text{d}t+ph(\tau)^{1/p}A^{p-1},
\end{align}
where $A:=M+\norm{f}_{\text{Ces}_p}+K^{1/p}$ and $K:=\sup_{s\in [0,1]}h(s)$.\par
From \eqref{opialcesaro eq:A5} and Fatou's Lemma it follows that 
\begin{align*}
&\lim_{n\to \infty}\norm{f_n}_{\text{Ces}_p(X)}^p \\
&\leq2^{p-1}\lim_{k\to \infty}\int_B\paren*{\frac{1}{t}\int_{[0,t]\cap C}\norm{f_{n_k}(s)-f(s)}\,\text{d}s}^p\,\text{d}t \\
&+2^{p-1}\lim_{k\to \infty}\int_B\paren*{\frac{1}{t}\int_{[0,t]\sm C}\norm{f_{n_k}(s)}\,\text{d}s}^p\,\text{d}t
+2^{p-1}\int_B\paren*{\frac{1}{t}\int_{[0,t]\cap C}g(s)\,\text{d}s}^p\,\text{d}t \\
&-2^{p-1}\int_B\paren*{\frac{1}{t}\int_{[0,t]\cap C}\liminf_{k\to \infty}\norm{f_{n_k}(s)-f(s)}\,\text{d}s}^p\,\text{d}t \\
&+\lim_{k\to \infty}\int_{[0,1]\sm B}\paren*{\frac{1}{t}\int_0^t\norm{f_{n_k}(s)}\,\text{d}s}^p\,\text{d}t.
\end{align*}
Combining this with \eqref{opialcesaro eq:A6} and \eqref{opialcesaro eq:A7} we obtain (by using $x^p+y^p\leq(x+y)^p$ for $x,y\geq 0$)
\begin{align*}
&\lim_{n\to \infty}\norm{f_n}_{\text{Ces}_p(X)}^p \\
&\leq2^{p-1}\lim_{k\to \infty}\norm{f_{n_k}-f}_{\text{Ces}_p(X)}^p+2^{p-1}ph(\tau)^{1/p}A^{p-1}+pL^{p-1}\eps^{1/p} \\
&+2^{p-1}\int_B\paren*{\frac{1}{t}\int_{[0,t]\cap C}g(s)\,\text{d}s}^p\,\text{d}t \\
&-2^{p-1}\int_B\paren*{\frac{1}{t}\int_{[0,t]\cap C}\liminf_{k\to \infty}\norm{f_{n_k}(s)-f(s)}\,\text{d}s}^p\,\text{d}t,
\end{align*}
thus
\begin{align}\label{opialcesaro eq:A8}
&\nonumber \lim_{n\to \infty}\norm{f_n}_{\text{Ces}_p(X)}^p \\
&\nonumber \leq2^{p-1}\lim_{n\to \infty}\norm{f_n-f}_{\text{Ces}_p(X)}^p+2^{p-1}ph(\tau)^{1/p}A^{p-1}+pL^{p-1}\eps^{1/p} \\
&+2^{p-1}\int_B\paren*{\paren*{\frac{1}{t}\int_{[0,t]\cap C}g(s)\,\text{d}s}^p-\paren*{\frac{1}{t}\int_{[0,t]\cap C}\varphi(s)\,\text{d}s}^p}\,\text{d}t.
\end{align}
Since $\lambda([0,1]\sm C)\leq\delta$ it follows from \eqref{opialcesaro eq:A3} that for $F\in \set*{g,\varphi}$ and $t\in (0,1-\frac{\tau}{3}]$ we have
\begin{equation*}
\abs*{\frac{1}{t}\int_0^tF(s)\,\text{d}s-\frac{1}{t}\int_{[0,t]\cap C}F(s)\,\text{d}s}\leq\frac{\theta}{t}
\end{equation*}
and hence
\begin{equation*}
\abs*{\paren*{\frac{1}{t}\int_0^tF(s)\,\text{d}s}^p-\paren*{\frac{1}{t}\int_{[0,t]\cap C}F(s)\,\text{d}s}^p}
\leq p\frac{\theta}{t}\paren*{\frac{1}{t}\int_0^tF(s)\,\text{d}s}^{p-1}.
\end{equation*}
Since $B\ssq [\frac{\tau}{3},1-\frac{\tau}{3}]$ it follows that
\begin{align*}
&\abs*{\int_B\paren*{\frac{1}{t}\int_0^tF(s)\,\text{d}s}^p\,\text{d}t-\int_B\paren*{\frac{1}{t}\int_{[0,t]\cap C}F(s)\,\text{d}s}^p\,\text{d}t} \\
&\leq\int_B p\frac{\theta}{t^p}\paren*{\int_0^tF(s)\,\text{d}s}^{p-1}\,\text{d}t\leq 
p\theta\paren*{\int_0^{1-\frac{\tau}{3}}F(s)\,\text{d}s}^{p-1}\int_{\frac{\tau}{3}}^{1-\frac{\tau}{3}}\frac{1}{t^p}\,\text{d}t \\
&=p\theta\paren*{\int_0^{1-\frac{\tau}{3}}F(s)\,\text{d}s}^{p-1}\frac{1}{1-p}\paren*{\paren*{1-\frac{\tau}{3}}^{1-p}-\paren*{\frac{\tau}{3}}^{1-p}}.
\end{align*}
Thus it follows from \eqref{opialcesaro eq:A2} that for $F\in \set*{g,\varphi}$ one has
\begin{equation}\label{opialcesaro eq:A9}
\abs*{\int_B\paren*{\frac{1}{t}\int_0^tF(s)\,\text{d}s}^p\,\text{d}t-\int_B\paren*{\frac{1}{t}\int_{[0,t]\cap C}F(s)\,\text{d}s}^p\,\text{d}t}\leq\eps.
\end{equation}
Since $\lambda([0,1]\sm B)\leq\tau$ we also have 
\begin{equation}\label{opialcesaro eq:A10}
\abs*{\int_B\paren*{\frac{1}{t}\int_0^tF(s)\,\text{d}s}^p\,\text{d}t-\int_0^1\paren*{\frac{1}{t}\int_0^tF(s)\,\text{d}s}^p\,\text{d}t}\leq\eps
\end{equation}
for $F\in \set*{g,\varphi}$, by \eqref{opialcesaro eq:A1}.\par
From \eqref{opialcesaro eq:A9} and \eqref{opialcesaro eq:A10} we obtain
\begin{equation*}
\abs*{\int_B\paren*{\frac{1}{t}\int_{[0,t]\cap C}F(s)\,\text{d}s}^p\,\text{d}t-\int_0^1\paren*{\frac{1}{t}\int_0^tF(s)\,\text{d}s}^p\,\text{d}t}\leq2\eps
\end{equation*}
for $F\in \set*{g,\varphi}$.\par
Together with \eqref{opialcesaro eq:A8} this implies
\begin{align*}
&\nonumber \lim_{n\to \infty}\norm{f_n}_{\text{Ces}_p(X)}^p \\
&\nonumber \leq2^{p-1}\lim_{n\to \infty}\norm{f_n-f}_{\text{Ces}_p(X)}^p+2^{p-1}ph(\tau)^{1/p}A^{p-1}+pL^{p-1}\eps^{1/p} \\
&+2^{p-1}\int_0^1\paren*{\paren*{\frac{1}{t}\int_0^tg(s)\,\text{d}s}^p-\paren*{\frac{1}{t}\int_0^t\varphi(s)\,\text{d}s}^p}\,\text{d}t+2^{p-1}4\eps.
\end{align*}
Hence by definition of $a$ we have
\begin{align*}
&\lim_{n\to \infty}\norm{f_n}_{\text{Ces}_p(X)}^p\leq2^{p-1}\lim_{n\to \infty}\norm{f_n-f}_{\text{Ces}_p(X)}^p \\
&+2^{p-1}ph(\tau)^{1/p}A^{p-1}+pL^{p-1}\eps^{1/p}-2^{p-1}a+2^{p+1}\eps.
\end{align*}
Since $h(\tau)\to 0$ for $\tau\to 0$ and $\tau<\eps$, we obtain for $\eps\to 0$
\begin{equation*}
\lim_{n\to \infty}\norm{f_n}_{\text{Ces}_p(X)}^p\leq2^{p-1}\lim_{n\to \infty}\norm{f_n-f}_{\text{Ces}_p(X)}^p-2^{p-1}a
\end{equation*}
and we are done.
\end{Proof}

Note that the assumptions that $X$ has the nonstrict Opial property and that $(f_n(t))_{n\in \N}$ converges weakly to zero a.\,e.
were only used to ensure that $a\geq 0$, which is only needed to conclude \eqref{eq:opialcesaro2} from \eqref{eq:opialcesaro1}. In other
words, \eqref{eq:opialcesaro1} is also valid without these two assumptions.\par
We have the following Corollary in the case that $X$ even has the Opial property (compare with \cite{hardtke}*{Corollary 4.2}).
\begin{corollary}\label{cor:opialcesaro}
Let $1\leq p<\infty$ and let $X$ be a Banach space with the Opial property. Let $(f_n)_{n\in \N}$ be a bounded sequence in $\mathrm{Ces}_p(X)$
such that $(f_n(t))_{n\in \N}$ converges weakly to zero for almost every $t\in [0,1]$. Suppose further that there exists a $g\in \mathrm{Ces}_p$ such
that $\norm{f_n(t)}\to g(t)$ a.\,e. Then 
\begin{equation*}
\limsup_{n\to \infty}\norm{f_n}_{\mathrm{Ces}_p(X)}<2^{1-1/p}\limsup_{n\to \infty}\norm{f_n-f}_{\mathrm{Ces}_p(X)} \ \ \forall f\in \mathrm{Ces}_p(X)\sm \set*{0}.
\end{equation*}
\end{corollary}

\begin{Proof}
Let $a$ be defined as in the previous proof. Since $X$ has the Opial property we have $\varphi(t)\geq g(t)$ for every $t\in [0,1]$ and even ``$>$'' if $f(t)\neq 0$, 
which by assumption happens on a set of positive measure. Thus $a>0$ and hence the desired inequality follows from Theorem \ref{thm:opialcesaro}.
\end{Proof}

Concerning the uniform Opial property, we also have the following analogue of \cite{hardtke}*{Theorem 4.3} for Ces\`aro function spaces
(the proof alos uses similar techbniques).
\begin{theorem}\label{thm:unifopialcesaro1}
Let $1\leq p<\infty$ and let $X$ be a Banach space with the uniform Opial property. Let $M,R>0$ and $f\in \mathrm{Ces}_p(X)\sm \set*{0}$.
Then there exists $\eta>0$ such that the following holds: whenever $(f_n)_{n\in \N}$ is a sequence in $\mathrm{Ces}_p(X)$ with 
$\sup_{n\in \N}\norm{f_n}_{\mathrm{Ces}_p(X)}\leq R$ such that $(f_n(t))_{n\in \N}$ converges weakly to zero and 
$\lim_{n\to \infty}\norm{f_n(t)}\leq M$ for almost every $t\in [0,1]$, then 
\begin{equation*}
\limsup_{n\to \infty}\norm{f_n}_{\mathrm{Ces}_p(X)}+\eta\leq2^{1-1/p}\limsup_{n\to \infty}\norm{f_n-f}_{\mathrm{Ces}_p(X)}.
\end{equation*}
\end{theorem}

\begin{Proof}
Fix $0<\tau<\norm{f}_{\text{Ces}_p(X)}$ and put $A:=\set*{s\in [0,1]:\norm{f(s)}\geq\tau}$. If $\lambda(A)=0$, then
we would obtain $\norm{f}_{\text{Ces}_p(X)}^p\leq\int_0^1 t^p\tau^p1/t^p\,\text{d}t=\tau^p$. Thus we must have $\lambda(A)>0$.
Let $w:=\eta_X(\tau,M)$.\par
Define $A_t:=A\cap [0,t]$ for $t\in [0,1]$. Then $\lambda(A_t)\to \lambda(A)$ for $t\to 1$ and hence we can find $t_0\in (0,1)$
such that $\lambda(A_t)\geq\lambda(A)/2$ for $t\in [t_0,1]$.\par
Put $\theta:=\int_{t_0}^11/t^p\,\text{d}t$ and $\nu:=\min\set*{(w^p\lambda(A)^p\theta/2)^{1/p},2^{1-1/p}(3R+1)}$.\par
Next we define $\omega:=2^{1-1/p}(3R+1)-(2^{p-1}(3R+1)^p-\nu^p)^{1/p}$ and finally $\eta:=\min\set*{\omega,1}$.\par
Now let $(f_n)_{n\in \N}$ be as above. Without loss of generality we may assume that $g(t):=\lim_{n\to \infty}\norm{f_n(t)}\leq M$ 
and $f_n(t)\to 0$ weakly for every $t\in [0,1]$. Let $\varphi(t):=\liminf\norm{f_n(t)-f(t)}$ for all $t\in [0,1]$. Then we have
$\varphi\geq g$ and the definition of $\eta_X$ implies that even $\varphi(s)-g(s)\geq\eta_X(\tau,M)=w$ for all $s\in A$.\par
Using the relation $(a-b)^p\leq a^p-b^p$ for $a\geq b\geq 0$ we obtain
\begin{align*}
&\paren*{\int_0^t\varphi(s)\,\text{d}s}^p-\paren*{\int_0^tg(s)\,\text{d}s}^p\geq\paren*{\int_0^t(\varphi(s)-g(s))\,\text{d}s}^p \\
&\geq\paren*{\int_{A_t}(\varphi(s)-g(s))\,\text{d}s}^p\geq w^p\lambda(A_t)^p
\end{align*}
for every $t\in [0,1]$.\par
Theorem \ref{thm:opialcesaro} now implies that 
\begin{align}\label{unifopialcesaro eq:1}
&\nonumber 2^{p-1}\limsup_{n\to \infty}\norm{f_n-f}_{\text{Ces}_p(X)}^p-\limsup_{n\to \infty}\norm{f_n}_{\text{Ces}_p(X)}^p \\
&\geq 2^{p-1}\int_0^1\frac{w^p}{t^p}\lambda(A_t)^p\,\text{d}t\geq2^{p-1}w^p\int_{t_0}^1\frac{\lambda(A_t)^p}{t^p}\,\text{d}t\geq w^p\frac{\lambda(A)^p}{2}\theta\geq\nu^p,
\end{align}
by the choice of $t_0$ and the definition of $\theta$ and $\nu$.\par
Next we define $h(s):=2^{1-1/p}s-(2^{p-1}s^p-\nu^p)^{1/p}$ for $s\geq2^{1/p-1}\nu$. It is easy to see that $h$ is decreasing on $[2^{1/p-1}\nu,\infty)$.\par
Now we consider two cases. If $\norm{f}_{\text{Ces}_p(X)}\geq 2R+1$, then 
\begin{align*}
&2^{1-1/p}\limsup_{n\to \infty}\norm{f_n-f}_{\text{Ces}_p(X)}\geq\limsup\norm{f_n-f}_{\text{Ces}_p(X)} \\
&\geq 2R+1-\liminf_{n\to \infty}\norm{f_n}_{\text{Ces}_p(X)}\geq 2R+1-R=R+1 \\
&\geq \limsup_{n\to \infty}\norm{f_n}_{\text{Ces}_p(X)}+\eta,
\end{align*}
since $\sup_{n\in \N}\norm{f_n}_{\mathrm{Ces}_p(X)}\leq R$ and $\eta\leq 1$.\par
If $\norm{f}_{\text{Ces}_p(X)}<2R+1$, then we have $\limsup\norm{f_n-f}_{\text{Ces}_p(X)}\leq 3R+1$ and hence (by \eqref{unifopialcesaro eq:1})
\begin{align*}
&\limsup_{n\to \infty}\norm{f_n}_{\text{Ces}_p(X)} \\
&\leq 2^{1-1/p}\limsup_{n\to \infty}\norm{f_n-f}_{\text{Ces}_p(X)}-h\paren*{\limsup\norm{f_n-f}_{\text{Ces}_p(X)}} \\
&\leq 2^{1-1/p}\limsup_{n\to \infty}\norm{f_n-f}_{\text{Ces}_p(X)}-h(3R+1) \\
&\leq 2^{1-1/p}\limsup_{n\to \infty}\norm{f_n-f}_{\text{Ces}_p(X)}-\eta,
\end{align*}
where the last inequality holds because $\eta\leq \omega=h(3R+1)$.
\end{Proof}

Finally, we have the following analogue of \cite{hardtke}*{Theorem 4.4} (we denote by $L^p([0,1],X)$ the $L^p$-Bochner space).
\begin{theorem}\label{thm:unifopialcesaro2}
Let $1<p<\infty$ and let $X$ be a Banach space with the uniform Opial property. Let $p<r\leq\infty$ and $\eps,M,K,R>0$. 
Then there exists $\eta>0$ such that the following holds: whenever $(f_n)_{n\in \N}$ is a sequence in $\mathrm{Ces}_p(X)$ 
with $\sup_{n\in \N}\norm{f_n}_{\mathrm{Ces}_p(X)}\leq R$ such that $(f_n(t))_{n\in \N}$ converges weakly to zero and 
$\lim_{n\to \infty}\norm{f_n(t)}\leq M$ for almost every $t\in [0,1]$ and $f\in L^r([0,1],X)\ssq L^p([0,1],X)\ssq \mathrm{Ces}_p(X)$ 
is such that $\norm{f}_r\leq K$ and $\norm{f}_{\mathrm{Ces}_p(X)}\geq\eps$, then
\begin{equation*}
\limsup_{n\to \infty}\norm{f_n}_{\mathrm{Ces}_p(X)}+\eta\leq2^{1-1/p}\limsup_{n\to \infty}\norm{f_n-f}_{\mathrm{Ces}_p(X)}.
\end{equation*}
\end{theorem}

\begin{Proof}
Let $s:=r/p\in(1,\infty]$ and let $s^{\prime}$ and $q$ be the conjugated exponents to $s$ and $p$. Choose $0<\tau<1$ such that
$q^p\tau^p<\eps^p$ and put $Q:=\min\set*{(\eps^p/q^p-\tau^p)^{s^{\prime}}K^{-ps^{\prime}},1}$, $w:=\eta_X(\tau,M)$ and $t_0:=1-Q/2$.\par
We also put $\theta:=\int_{t_0}^1 1/t^p\,\text{d}t$ and $\nu:=\min\set*{(w^pQ^p\theta/2)^{1/p},2^{1-1/p}(3R+1)}$,
as well as $\omega:=2^{1-1/p}(3R+1)-(2^{p-1}(3R+1)^p-\nu^p)^{1/p}$ and finally $\eta:=\min\set*{\omega,1}$.\par
Now let $(f_n)_{n\in \N}$ in $\text{Ces}_p(X)$ and $f\in L^r([0,1],X)$ be as above. We assume without loss of generality that 
$g(t):=\lim_{n\to \infty}\norm{f_n(t)}\leq M$ and $f_n(t)\to 0$ weakly for every $t\in [0,1]$.\par
Let $A:=\set*{s\in [0,1]:\norm{f(s)}\geq\tau}$. Since $\eps\leq\norm{f}_{\text{Ces}_p(X)}\leq q\norm{f}_p$ (see (4) on page 2)
we can proceed analogously to the proof of \cite{hardtke}*{Theorem 4.4} to show that $\lambda(A)\geq Q$.\par
Let $A_t:=A\cap [0,t]$ for $t\in [0,1]$. We have $\lambda(A)-\lambda(A_{t_0})=\lambda(A\cap (t_0,1])\leq 1-t_0=Q/2$ and hence 
$\lambda(A_t)\geq \lambda(A_{t_0})\geq Q/2$ for $t\in [t_0,1]$.\par
As in the previous proof we can now use Theorem \ref{thm:opialcesaro} to conclude
\begin{equation*}
2^{p-1}\limsup_{n\to \infty}\norm{f_n-f}_{\text{Ces}_p(X)}^p-\limsup_{n\to \infty}\norm{f_n}_{\text{Ces}_p(X)}^p\geq\nu^p
\end{equation*}
and from this obtain, also as in the previous proof, that 
\begin{equation*}
\limsup_{n\to \infty}\norm{f_n}_{\text{Ces}_p(X)}\leq 2^{1-1/p}\limsup_{n\to \infty}\norm{f_n-f}_{\text{Ces}_p(X)}-\eta.
\end{equation*}
\end{Proof}

\address
\email

\end{document}